\documentclass{amsart}

\usepackage{amsmath,amsthm,amssymb,amscd,enumerate,mathtools,enumitem,tikz-cd,bm}
\usepackage{amsfonts}
\usepackage{rotating}
\usepackage{euscript}
\usepackage{pst-node}
\usepackage{epsfig,verbatim}

\def\be{\begin{equation}}
\def\ee{\end{equation}}

\def\C{{\mathbb C}}

\def\P{{\mathbb P}}

\def\phi{{\varphi}}
\def\v{{\varepsilon}} 

\def\deg{{\rm deg\,}}

\def\GCD{{\rm GCD }}

\def\bp{\begin{proposition}}
\def\ep{\end{proposition}}

\def\bt{\begin{theorem}}
\def\et{\end{theorem}}
\def\br{\begin{remark}}
\def\er{\end{remark}}
\def\be{\begin{equation}}
\def\bee{\begin{equation*}}
\def\l{\label}
\def\la{\label}

\def\ee{\end{equation}}
\def\eee{\end{equation*}}
\def\bl{\begin{lemma}}
\def\el{\end{lemma}}
\def\bc{\begin{corollary}}
\def\ec{\end{corollary}}
\def\pr{\noindent{\it Proof. }}

\def\bd{\begin{definition}}
\def\ed{\end{definition}}
\def\t{\widetilde}
\def\h{\widehat}
\def\hat{\widehat}
\def\tilde{\widetilde}

\def\t{\widetilde }

\newtheorem{theorem}{Theorem}[section]
\newtheorem{lemma}[theorem]{Lemma}
\newtheorem{definition}[theorem]{Definition}
\newtheorem{corollary}[theorem]{Corollary}
\newtheorem{proposition}[theorem]{Proposition}
\newtheorem{problem}[theorem]{Problem}

\theoremstyle{definition}
\newtheorem{example}[theorem]{Example}
\theoremstyle{definition}
\newtheorem{remark}[theorem]{Remark}
\def\bpr{\begin{problem}}
\def\epr{\end{problem}}

\mathtoolsset{showonlyrefs}

\begin{document}

\title{On intertwined polynomials}

\author[F. Pakovich]{Fedor Pakovich}
\thanks{
This research was supported by ISF Grant No. 1092/22} 
\address{Department of Mathematics, 
Ben Gurion University of the Negev, P.O.B. 653, Beer Sheva,  8410501, Israel}
\email{pakovich@math.bgu.ac.il}

\begin{abstract}
Let $A_1$ and $A_2$ be polynomials of degree at least two over $\C$. We say that $A_1$ and $A_2$ are intertwined if the endomorphism $(A_1, A_2)$ of $\C\P^1 \times \C\P^1$ given by  $(z_1, z_2) \mapsto (A_1(z_1), A_2(z_2))$ admits an irreducible periodic curve that is neither a vertical nor a horizontal line. For a polynomial $A$, we denote by $\mathrm{Inter}(A)$ the set
of all polynomials $B$ such that some iterate of $B$ is intertwined
with some iterate of $A$. 
In this paper, we prove a conjecture of Favre and Gauthier describing the structure of $\mathrm{Inter}(A)$. We also obtain a sharp bound on the periods of periodic curves for endomorphisms $(A_1, A_2)$ in terms of the sizes of the symmetry groups of the Julia sets of $A_1$ and $A_2$.

\end{abstract}

\maketitle

\section{Introduction} 
In their seminal paper, Medvedev and Scanlon \cite{ms} described  invariant curves of endomorphisms of 
$(\C\P^1)^2$ of the form
\be \l{end} 
(A_1, A_2)\colon (z_1, z_2) \mapsto (A_1(z_1), A_2(z_2)),
\ee
where $A_1$ and $A_2$ are polynomials, and showed that describing invariant curves for similar endomorphisms of $(\C\P^1)^m$ with $m>2$ reduces to the case $m=2$.  
Since then, these results have been widely used in various contexts of 
complex and arithmetic dynamics (see, e.g., 
\cite{bm}, \cite{br}, \cite{fg}, \cite{favg}, \cite{gn}, \cite{gn1}, \cite{gn2}, 
\cite{gx}, \cite{n}).

From an algebraic perspective, the description of invariant—and
more generally, periodic—curves for endomorphisms \eqref{end} essentially 
reduces to a system of polynomial semiconjugacies
involving iterates of $A_1$ and $A_2$.
In \cite{ms}, such systems were analyzed using Ritt's theory
of polynomial decompositions \cite{r1}. An alternative approach, combining Ritt's theory with results from \cite{p1}, was proposed in \cite{pj}. 
In a sense, this approach is less technical and also allows one, to some extent, to complete the results of \cite{ms}. 
Nevertheless, the approaches of both \cite{ms} and \cite{pj} do not extend to the case where $A_1$ and $A_2$ are arbitrary  rational functions. In this setting, a description  of invariant and periodic curves for endomorphisms \eqref{end} was obtained in \cite{ic}, based on a characterization of semiconjugacies between rational functions provided in a series of papers \cite{semi}, \cite{dyna}, \cite{rec}, \cite{lattes}, \cite{fin}.

The goal of this paper remains the study of periodic curves for polynomial endomorphisms \eqref{end}, but from a more general perspective. Here, we fix a polynomial $A_1$ and aim to describe all polynomials $A_2$ such that, for some {\it iterates} of $A_1$ and $A_2$, the corresponding endomorphism \eqref{end} admits an irreducible periodic curve. 
More precisely, we say that polynomials $A_1$ and $A_2$ of degree at least two are \emph{intertwined} if the endomorphism \eqref{end} admits an irreducible periodic curve that is neither a vertical nor a horizontal line. Given a polynomial $A$, we denote by $\mathrm{Inter}(A)$ the collection of all polynomials $B$ for which there exist $k, l \ge 1$ such that $A^{\circ k}$ and $B^{\circ l}$ are intertwined. 

We denote by $J(A)$ the Julia set of $A$, and by 
$\Sigma(A)$ the symmetry group of $J(A)$, consisting of all affine transformations $\mu = az + b$ that leave $J(A)$ invariant. 
Note that the group $\Sigma(A)$ is finite unless $A$ is conjugate to a power. 

In the above notation, our main result is the following statement, conjectured by Favre and Gauthier \cite[p.~88]{favg}.

\bt \l{t2}
For any $d \ge 2$, there exists a constant $r(d)$ such that for every polynomial $A$ of degree $d$, there are polynomials $B_1, \dots, B_r$ with $r \le r(d)$ satisfying the following property: if $B \in \mathrm{Inter}(A)$, then $B$ is conjugate to $\mu \circ B_i^{\circ n}$ for some $i$, $1 \le i \le r$, some $\mu \in \Sigma(B_i)$, and  some  $n \ge 1$.
\et

Note that \cite{favg} uses a different definition of $\Sigma(A)$,
which yields the same group as ours unless $A$ is conjugate
to a power (see \cite[Definition~3.1 and Proposition~3.9]{favg}).
For $A=z^d$, the group in \cite{favg} consists of the rotations $z\mapsto\zeta z$
for roots of unity $\zeta$ and differs from ours.
Nevertheless, Theorem~\ref{t2} remains valid under either definition.

Our proof of Theorem \ref{t2} yields an effective bound for $r(d)$.
However, the resulting estimate does not appear to be
optimal, so we do not make it explicit.

As a by-product of our method, we obtain a bound on the period of an  irreducible periodic curve with nonconstant coordinate projections for the endomorphisms \eqref{end} in the case when $A_1$ and $A_2$ are not special—that is, not conjugate to $z^d$ or to $\pm T_d$—in terms of the sizes of $\Sigma(A_1)$ and $\Sigma(A_2)$. 
If such a curve exists, then the degrees of $A_1$ and $A_2$ 
coincide, and the explicit bound $d^{2d^4}/2$ 
 in terms of this common degree $d$ 
was established in \cite[Theorem~1.1, p.~214]{gn1}.

We prove the following result.

\bt \l{t1}
Let $A_1$ and $A_2$ be non-special polynomials of degree $d\geq2$,
and let $C$ be an irreducible $(A_1,A_2)$-periodic curve
that is neither a vertical nor a horizontal line.
Then its period is at most
$
\operatorname{lcm}\bigl(|\Sigma(A_1)|,|\Sigma(A_2)|\bigr).
$ 
\et

Note that the bound of Theorem~\ref{t1}  is attained for polynomials of arbitrarily
large degree with arbitrarily large symmetry groups
(see Example~\ref{sharpperiod}).

Since $|\Sigma(A_i)|\leq d$, Theorem~\ref{t1} gives the bound
$d^2$ in terms of the common degree, improving on the bound
obtained in \cite{gn1}. It also yields the following interesting corollary
for polynomials with trivial symmetry groups
(note that such polynomials necessarily have degree
greater than two).

\bc Let $A_1$ and $A_2$ be non-special polynomials of degree $d> 2$ such that  the groups $\Sigma(A_1)$ and $\Sigma(A_2)$ are trivial. Then any irreducible $(A_1,A_2)$-periodic curve
that is neither a vertical nor a horizontal line is invariant.
\ec

\smallskip

In this paper, we work over $\mathbb C$, and our proofs use
complex analytic methods. The situation in positive
characteristic involves additional phenomena
(see, e.g., \cite[Remark~2.31]{ms}) and requires separate analysis,
which we do not undertake here. 

The paper is organized as follows. In Section~2, we first recall two 
classical results concerning the functional equation
\[
A \circ C = D \circ B,
\]
where \(A\), \(B\), \(C\), and \(D\) are polynomials. We then review basic
results on the functional equation
\[
A \circ X = X \circ B,
\]
which describes the semiconjugacy relation for polynomials, providing sketch
proofs throughout.

The material of the next section is also known, but we provide complete proofs, 
as we believe the way these results are presented may be of independent interest. 
We begin by proving a particular case of the general result stated in \cite{p1}. 
This case characterizes polynomials $A_1$, $A_2$ and infinite compact sets 
$K_1, K_2 \subset \C$ satisfying  
\be \l{des} 
A_1^{-1}(K_1) = A_2^{-1}(K_2),
\ee  
under the condition $\deg A_1 \mid \deg A_2$. From this result, we deduce a 
criterion for semiconjugacy of polynomials in terms of their Julia sets, first proved in \cite{pj}. 
We also recover several  results concerning the group $\Sigma(A)$ and 
polynomials sharing Julia sets, proved in \cite{b}, \cite{be2}, \cite{be1}, 
and \cite{sh}.

In Section~4 we establish several results concerning the implications
of the semiconjugacy relation for polynomials in the case when one of
the polynomials involved is an $l$th iterate of some polynomial. In
Section~5, building on these findings, we prove Theorem~\ref{t1}. In
Section~6 we prove additional results about solutions to 
 semiconjugacy relations between iterates $A^{\circ l}$,
$l \ge 1$, and $B^{\circ s}$, $s \ge 1$, assuming that one of $A$ or
$B$ is fixed while the other, together with the integers $s$ and $l$,
may vary. Finally, we prove Theorem~\ref{t2}.

\section{Polynomial semiconjugacies} 
\subsection{Theorems of Engstrom and Ritt}
We begin by recalling two results about polynomial solutions to the equation \be \la{ura} A \circ C = D \circ B. \ee

The first result, proved by Engstrom  \cite{en}, roughly speaking, reduces the problem of finding solutions to \eqref{ura} to the case \be \l{ggccdd} \text{GCD}(\deg A, \deg D) = 1, \quad \text{GCD}(\deg C, \deg B) = 1. \ee

\bt\la{r1}
Let $A,C,D,B$ be non-constant polynomials 
such that \eqref{ura} holds. Then there exist polynomials
$U, V, \widetilde A, \widetilde C, \widetilde D, \widetilde B,$ where
$$\deg U=\GCD(\deg A,\deg D),  \ \ \ \deg V=\GCD(\deg C,\deg B),$$
such that
$$A=U\circ \widetilde A, \ \  D=U\circ \widetilde D, \ \ C=\widetilde C\circ V, \ \  B=\widetilde B\circ V,$$
and 
$$ \widetilde A\circ \widetilde C=\widetilde D\circ \widetilde B.$$ \et 

Notice that Theorem~\ref{r1} implies that if $\deg C = \deg B$ in \eqref{ura}, then there exists a degree-one polynomial $\mu$ such that 
\[
A = D \circ \mu^{-1}, \qquad C = \mu \circ B.
\]

The second result, proved by Ritt \cite{r1}, describes polynomial solutions to \eqref{ura} that satisfy \eqref{ggccdd}.
 
\bt \la{r2}
Let $A,C,D,B$ be non-constant polynomials  
such that \eqref{ura} and \eqref{ggccdd} hold.
Then there exist polynomials $\sigma_1,\sigma_2,\mu, \nu$ of degree one  
such that, up to a possible replacement of $A$ by $D$ and of $C$ by $B$, either
\begin{align}  \la{rp11}  &A=\nu \circ z^sR^n(z) \circ \sigma_1^{-1}, & &C=
\sigma_1 \circ z^n \circ \mu \\ 
&\la{rp11+} D=\nu \circ z^n \circ \sigma_2^{-1},& 
&B=\sigma_2 \circ
z^sR(z^n) \circ \mu, \end{align}
where $R$ is a polynomial, $n\geq 1,$ $s\geq 0$,  and $\GCD(s,n)=1,$ or \pagebreak 
\begin{align} \la{rp21} &A=\nu \circ T_m \circ \sigma_1^{-1}, &
&C=\sigma_1 \circ T_n \circ \mu, \\
&\la{rp21+} D=\nu \circ T_n \circ \sigma_2^{-1},& &B=\sigma_2 \circ
T_m\circ \mu, \end{align} 
where $T_n, T_m$ are the Chebyshev polynomials, $n,m\geq 1$, and $\GCD(n,m)=1.$ 
\et

Notice that although formulas \eqref{rp11} and \eqref{rp21} hold up to a possible replacement of $A$ by $D$ and of $C$ by $B$, we may assume the choice given by \eqref{rp11} and \eqref{rp21} whenever $m = \deg A = \deg B$ is greater than $n = \deg C = \deg D$. Indeed, if $R$ is not a constant, this is obvious. On the other hand, if $R$ is a constant, then by adjusting $\nu$ and $\sigma_2$, we can assume that $R = 1$, and formulas \eqref{rp11} become symmetric, like formulas \eqref{rp21}.

Notice also  that if $m > n > 2$, then \eqref{rp11} and \eqref{rp21} cannot hold simultaneously since the multiplicity of $T_n$ at any point of $\mathbb{C}$ is at most  two. 
However, as the following example shows, this is not the case for $n = 2$, which sometimes makes the application of Theorem~\ref{r2} difficult.

\begin{example} 
The formula 
\be \l{equ}
T_n = \frac{n}{2} \sum_{k = 0}^{\lfloor n/2 \rfloor}
\frac{(-1)^k (n-k-1)!}{k!\,(n-2k)!} (2x)^{n-2k} 
\ee
implies that for odd $n$ the equality
\[
T_n(z) = z E_n(z^2) 
\]
holds for some polynomial $E_n$. Furthermore, since $T_2 = \theta \circ z^2$, where $\theta = 2z - 1$, we have
\[
z E_n(z^2) \circ \theta \circ z^2 = T_n \circ T_2 = T_2 \circ T_n = \theta \circ T_n^{2} =
\theta \circ z^2 E_n^2(z^2) 
=\theta \circ z E_n^2(z) \circ z^2.
\]
Since the last equality implies
\[
z E_n(z^2) \circ \theta = \theta \circ z E_n^2(z),
\]
we conclude that
\[
T_n = \theta \circ z E_n^2(z) \circ \theta^{-1}. 
\]
Therefore, the equality $$T_n \circ T_2 = T_2 \circ T_n$$ may be written in the form
\[
(\theta \circ z E_n^2(z) \circ \theta^{-1}) \circ (\theta \circ z^2) = (\theta \circ z^2) \circ z E_n(z^2).
\]
\end{example}

\subsection{Primitive semiconjugacies} 

Let \( A \) and \( B \) be polynomials of degree at least two.  
Let us recall that if the equality  
\begin{equation} \label{1}
A \circ X = X \circ B
\end{equation}   
is satisfied for some polynomial \( X \) of degree one, the polynomials \( A \) and \( B \) are called {\it conjugate}. If \eqref{1} is satisfied for some polynomial \( X \) of degree at least two, the polynomial \( B \) is called {\it semiconjugate} to \( A \).  

We will use the notation \( A \leq B \) if, for polynomials \( A \) and \( B \), there exists a non-constant polynomial \( X \) such that \eqref{1} holds, and the notation \( A \underset{X}{\leq} B \) if \( A \), \( B \), and \( X \) satisfy \eqref{1}.  
Notice that by writing equality \eqref{1} in the form of a commuting diagram 
\be \l{11} 
\begin{CD} 
\C\P^1 @> B >> \C\P^1 \\ 
@V X VV @VV X V\\ 
\C\P^1 @> A >> \C\P^1,
\end{CD} 
\ee
we immediately see that \( A \underset{X}{\leq} B \) implies \( A^{\circ i} \underset{X}{\leq} B^{\circ i} \) for every \( i \geq 2 \).

Theorem \ref{r2} allows us to describe polynomial solutions of \eqref{1}   under the condition  
\be \l{asdf}
\text{GCD}(\deg X, \deg B) = 1.
\ee  
We call such solutions of \eqref{1} {\it primitive}. Notice that Theorem \ref{r1}, 
together with L\"uroth's theorem, implies that a solution $A, B, X$ of \eqref{1} 
is primitive if and only if  
\be \l{con}
\C(X, B) = \C(z).
\ee  
Below, we will use the conditions \eqref{asdf} and \eqref{con} interchangeably.

Specifically, Theorem \ref{r2} implies the following result  (cf. \cite{i}).

\bt \la{i}  
Let $A, B, X$ be polynomials of degree at least two such that \( A \underset{X}{\leq} B \) and $\GCD(\deg X,\deg B)=1$. Then there exist polynomials $\mu, \nu$ of degree one  such that either  
\be \l{ifab} A = \nu \circ z^s R^n(z) \circ \nu^{-1}, \quad X = \nu \circ z^n \circ \mu, \quad B = \mu^{-1} \circ z^s R(z^n) \circ \mu,\ee  
where $R$ is a polynomial, $n \ge 2$, $s \ge 1$, and $\gcd(s,n) = 1$, or  
\be \l{ifab2} A = \nu \circ (\pm 1)^n T_m \circ \nu^{-1}, \quad X = \nu \circ T_n \circ \mu, \quad B = \mu^{-1} \circ \pm T_m \circ \mu,\ee  
where $T_n, T_m$ are Chebyshev polynomials, $n, m \ge 2$, and $\gcd(n, m) = 1.$  
\et

\pr Applying Theorem~\ref{r2} to the relation
$
\t A \underset{X}{\leq} \t B, 
$
where 
\be \l{cre} \t A = A^{\circ i}, \qquad \t B = B^{\circ i},\ee 
for an integer $i\geq1$ such that $(\deg B)^i>\deg X$, we conclude that, without loss of generality, 
we may assume that either $X = T_n$ with $n \geq 3$, or $X = z^n$ with $n \geq 2$. 
Moreover, in the first case, there exist degree-one polynomials $\nu$ and $\sigma$ such that 
\be \l{ci2}
\t A = \nu \circ T_{m^i}, \qquad \t B = \sigma \circ T_{m^i}, 
\qquad
X = T_n = \nu \circ T_n \circ \sigma^{-1}.
\ee

It follows from formula \eqref{equ} that \(T_n\) has nonzero coefficients \(c_n\) and \(c_{n-2}\) in degrees \(n\) and \(n-2\), while the coefficient of \(z^{n-1}\) vanishes. Thus, if the equality for $X$ in \eqref{ci2} holds, then writing $$\sigma^{-1} = az + b\quad  {\rm and}  \quad \nu = cz + d,$$ we see that $b = 0$.
Furthermore, since $n-2 > 0$, we have  
\[
c c_n a^n = c_n \quad \text{and} \quad c c_{n-2} a^{n-2} = c_{n-2},
\]  
which implies $a = \pm 1$. Thus, $\sigma^{-1} = \sigma = \pm z$, and consequently $\nu = \pm z$.  
Since polynomials of the form $\pm T_m$ with $m \geq 2$ are uniquely determined by the property that their Julia set is the segment $[-1,1]$, it follows now from the equalities \eqref{cre} that $A = \pm T_m$ and $B = \pm T_m$. Finally, \eqref{1} determines the choice of the sign
in \eqref{ifab2}.

In the case where $X = z^n$ with $n \geq 2$, we finish the proof using Theorem~\ref{r1} as follows. It follows from \eqref{1} that 
\[
z^n \circ B = z^n \circ (B \circ \v_n z),
\] 
where $\v_n$ denotes a primitive $n$th root of unity. Applying Theorem~\ref{r1} to this equality, we conclude that there exists a degree-one polynomial $\mu$ such that the equalities 
\[
z^n=z^n \circ \mu^{-1}, \qquad   B \circ \v_n z = \mu \circ B
\] hold. 
The first equality implies that $$\mu = \v_n^s z$$ for some $s \geq 0$. The second equality then implies that $B=z^s R(z^n)$
for some polynomial $R$, and \eqref{1} gives $A=z^s R^n(z)$.
Moreover, $\gcd(s,n)=1$ because $\gcd(n,m)=1$.
 \qed

Notice that Theorem~\ref{i} is a special case of the description of solutions to \eqref{1} in arbitrary rational functions satisfying condition \eqref{con}, as obtained in \cite{semi}. For further details, we refer the reader to \cite{semi}, \cite{dyna}, and \cite{lattes}. An important distinction of this general setting is that conditions \eqref{asdf} and \eqref{con} are no longer equivalent, since Theorem~\ref{r1} does not hold. 

\begin{example}\l{ex1}
The simplest counterexample to Theorem~\ref{r1} in the rational case is given by the equality  
\be \l{rr1} 
\left( \frac{z^2 - 1}{z^2 + 1} \right) \circ 
\left( \frac{z^2 - 1}{z^2 + 1} \right)
= 
- \frac{2}{z^2 - 2} \circ \left( z + \frac{1}{z} \right).
\ee  
Here all functions involved have degree two, yet the equality  
\be \l{rr2} 
\frac{z^2 - 1}{z^2 + 1} = \mu \circ \left( z + \frac{1}{z} \right)
\ee  
for any M\"obius transformation $\mu$ is impossible, because the two sides have 
different sets of critical points.
\end{example}

\subsection{Elementary transformations}

In this section, we recall how an arbitrary solution of \eqref{1} can be reduced to a primitive one. 
Let $A$ be a polynomial. A polynomial $B$ is called an \emph{elementary transformation} of $A$ if there exist polynomials $U$ and $V$ such that
\[
A = U \circ V \quad \text{and} \quad B = V \circ U.
\]
Clearly, the identity
\be \label{ii}
(U \circ V) \circ U = U \circ (V \circ U)
\ee
implies that \( A \underset{U}{\leq} B \). Similarly, \( B \underset{V}{\leq} A \).

We say that polynomials $A$ and $B$ are \emph{equivalent}, and write $A \sim B$, if there exists a chain of elementary transformations connecting $B$ and $A$. If
\be \label{tr}
B \to B_1 \to B_2 \to \cdots \to B_s = A
\ee
is such a chain, and if \( U_i, V_i \), \( 1 \le i \le s \), are the corresponding polynomials satisfying
\be \label{pol}
B = V_1 \circ U_1, \quad B_i = U_i \circ V_i,\quad 1 \le i \le s,
\ee
and
\[
U_i \circ V_i = V_{i+1} \circ U_{i+1},\quad 1 \le i \le s - 1,
\]
then the composition
\[
X = U_s \circ U_{s-1} \circ \cdots \circ U_1
\]
makes the diagram \eqref{11}  
commutative. Thus, \( A \sim B \) implies both \( A \leq B \) and \( B \leq A \).

Notice that, for any polynomial $X$ of degree one, the equality  
$$A = (A \circ X) \circ X^{-1}$$  
implies  
$$A \sim X^{-1} \circ A \circ X,$$  
so that each equivalence class is a union of conjugacy classes. Moreover, the number of such 
classes is finite and can be bounded solely in terms of the degree of $A$ (see \cite{rec}, 
\cite{fin} for more details).  

Notice also that since any decomposition $V \circ U$ of $B$ naturally induces a decomposition $(V \circ U)^{\circ e}$ of $B^{\circ e}$ for $e \geq 2$, and since \eqref{ii} implies 
\[
(U \circ V)^{\circ e} \circ U = U \circ (V \circ U)^{\circ e},
\]
the relation $A \sim B$ implies the relation $A^{\circ e} \sim B^{\circ e}$ for every $e \geq 2$.

An arbitrary solution of equation \eqref{1} can be reduced to a primitive one by a sequence of elementary transformations. Specifically, if \( \mathbb{C}(X, B) \neq \mathbb{C}(z) \), then the Lüroth theorem implies that there exists a polynomial \( W \) of degree greater than one such that
\be \label{of}
B = \widetilde{B} \circ W, \quad X = \widetilde{X} \circ W
\ee
for some polynomials \( \widetilde{X} \) and \( \widetilde{B} \) with \( \mathbb{C}(\widetilde{X}, \widetilde{B}) = \mathbb{C}(z) \). Moreover, the diagram
\[
\begin{CD}
\mathbb{C}\mathbb{P}^1 @> B >> \mathbb{C}\mathbb{P}^1 \\
@V W VV @VV W V \\
\mathbb{C}\mathbb{P}^1 @> W \circ \widetilde{B} >> \mathbb{C}\mathbb{P}^1 \\
@V \widetilde{X} VV @VV \widetilde{X} V \\
\mathbb{C}\mathbb{P}^1 @> A >> \mathbb{C}\mathbb{P}^1
\end{CD}
\]
commutes. Hence, the triple \( A, \widetilde{X}, W \circ \widetilde{B} \) is another solution of \eqref{1}. While this new solution may still be non-primitive, it satisfies $$ \deg \widetilde{X} < \deg X.$$ Therefore, by repeating this procedure, we eventually obtain a primitive solution. 

Thus, we conclude the following.

\begin{lemma} \l{25}
Let \(A\) and \(B\) be polynomials of degree at least two, and let \(X\) be a non-constant polynomial such that \( A \underset{X}{\leq} B \). Then there exist polynomials \( X_0 \), \( B_0 \), and \( U \) such that:
\begin{enumerate}[label=\upshape(\alph*)]
    \item The diagram
    \[
    \begin{CD}
    \mathbb{C}\mathbb{P}^1 @> B >> \mathbb{C}\mathbb{P}^1 \\
    @V U VV @VV U V \\
    \mathbb{C}\mathbb{P}^1 @> B_0 >> \mathbb{C}\mathbb{P}^1 \\
    @V X_0 VV @VV X_0 V \\
    \mathbb{C}\mathbb{P}^1 @> A >> \mathbb{C}\mathbb{P}^1
    \end{CD}
    \]
    commutes, and \( X = X_0 \circ U \),

    \item The relation \( B_0 \sim B \) holds,

    \item The triple \( A, B_0, X_0 \) is a primitive solution of equation~\eqref{1}. \qed
\end{enumerate}
\end{lemma}

Notice that for arbitrary  rational functions $A$, $B$, and $X$,
Lemma \ref{25} remains valid if one defines a primitive solution of
\eqref{1} as a solution satisfying \eqref{con} and extends the
definition of an elementary transformation by allowing $U$ and $V$ to
be rational as well (for more details, see \cite{rec}, \cite{lattes}).

\subsection{\l{24} Semiconjugacies between special polynomials} 
We recall that a polynomial \(A\) of degree at least two is called \emph{special} if it is conjugate either to \(z^d\) or to \(\pm T_d\). It is well known that special polynomials are characterized by their Julia sets: \(A\) is conjugate to \(z^d\) if and only if   $J(A)$ is a circle, and \(A\) is conjugate to \(\pm T_d\) if and only if  $J(A)$ is a segment. In particular, \(\pm T_d\) is not conjugate to \(z^d\). Moreover, \(A\) is conjugate to \(\pm T_d\) (respectively \(z^d\)) if and only if some iterate of \(A\) is.

The following result was established in \cite{pj}.
Below, we give an alternative proof using Lemma~\ref{25} and Theorem~\ref{i}.

\bt \l{spe} 
Let $A$ and $B$ be polynomials of degree at least two such that $A \leq B$. Then $A$ is conjugate to $z^m$ if and only if $B$ is conjugate to $z^m$. 
Similarly, $A$ is conjugate to $\pm T_m$ if and only if $B$ is conjugate to $\pm T_m$.
\et 

\pr 
It  follows  from Theorem~\ref{r1} that if 
$T_m = U_1 \circ U_2$ for some polynomials $U_1$ and $U_2$ of degrees $m_1$ and $m_2$, respectively, 
then there exists a polynomial $\alpha$ of degree one such that 
\[
U_1 = T_{m_1} \circ \alpha, \qquad U_2 = \alpha^{-1} \circ T_{m_2}.
\]
Similarly, the equality $z^m = U_1 \circ U_2$ implies that 
\[
U_1 = z^{m_1} \circ \alpha, \qquad U_2 = \alpha^{-1} \circ z^{m_2}.
\]  
Therefore, a polynomial is conjugate to $z^m$ or $\pm T_m$ if and only if all polynomials in its equivalence class share this property. 

Combining this with Lemma~\ref{25}, we see that it suffices to prove the theorem under the assumption that there exists a polynomial $X$ 
 with $\deg X \geq 2$ such that 
$A, B,$ and $ X$ form a primitive solution of equation~\eqref{1}. Moreover, 
 we may assume $\deg X > 2$. Indeed, choose an odd integer $l>1$ with $\gcd(l,m)=1$.
Since $\pm T_m$ commutes with $T_l$, the relations 
$$\pm T_m \underset{X}{\leq} B\quad \text{and} \quad A \underset{X}{\leq} \pm T_m$$ imply the relations 
\[
\pm T_m \underset{T_l \circ X}{\leq} B \quad \text{and} \quad A \underset{X \circ T_l}{\leq} \pm T_m,
\]
and similar formulas hold for $z^m$ and $z^l$ in place of $\pm T_m$ and $T_l$.

Since $z^m$ and $\pm T_m$ are not conjugate, it follows from Theorem~\ref{i} that, to prove the theorem, it is enough to show that, in the formulas \eqref{ifab}, $A$ is conjugate to $z^m$ if and only if $B$ is conjugate to $z^m$, and that, in the case of non‑constant $R$, neither $A$ nor $B$ can be conjugate to $\pm T_m$. 

The first statement is straightforward.
To prove the second, observe that if $A$ is conjugate to
$\pm T_m$, then $n=2$, since the multiplicity of $T_m$
at any point of $\mathbb{C}$ is at most two.
This contradicts $\deg X>2$.
On the other hand, assuming, after an affine change of
coordinates, that $B=\pm T_m$, we obtain from \eqref{equ}
that $\mu(0)=0$ and $n=2$, again contradicting $\deg X>2$.
\qed

\vskip 0.2cm

Note that yet another proof of Theorem \ref{spe} can be obtained by applying
\cite[Theorem~1.14]{jx}, using the identity
\[
(A^{\circ n})'(X(p))=(B^{\circ n})'(p),
\]
which follows from \eqref{1} whenever $B^{\circ n}(p)=p$
and $X$ is unramified at $p$.
To distinguish the two special types, observe that polynomial
semiconjugacy preserves the existence of a finite superattracting
periodic point, which exists for $z^m$ but not for $\pm T_m$.

The same theorem also implies that, if \eqref{1} holds
with a nonconstant {\it rational} function $X$, then $A$ is
special if and only if $B$ is special, although one may
be conjugate to $z^m$ and the other to $\pm T_m$.

\section{Semiconjugacies and Julia sets} 
\subsection{Polynomials sharing preimages of compact sets} 
Our approach to polynomial semiconjugacies is based on the characterization provided by Theorem~\ref{j+} below, first established in \cite{pj}. 
This result is a corollary of the results of \cite{p1}, which describe polynomials \(A_1, A_2\) and compact sets \(K_1, K_2\) satisfying
\be \l{eqv}
A_1^{-1}(K_1) \,=\, A_2^{-1}(K_2),
\ee
in terms of solutions of the functional equation
\be \l{2}
B_1 \circ A_1 \,=\, B_2 \circ A_2,
\ee
where \(A_1, A_2, B_1, B_2\) are polynomials.

For any polynomial solution of \eqref{2} and any compact set 
\(K \subset \C\), one obtains a solution of \eqref{eqv} by setting  
\be \l{3}
K_1 = B_1^{-1}(K), \qquad K_2 = B_2^{-1}(K),
\ee
and, roughly speaking, the main result of \cite{p1} states that all solutions 
of \eqref{eqv} arise in this way whenever the set defined by \eqref{eqv} contains 
at least \(\operatorname{LCM}(\deg A_1, \deg A_2)\) points.

For the sake of completeness, we provide a full proof of the result of \cite{p1} in the 
special case where $\deg A_1$ divides $\deg A_2$ (Theorem \ref{be}), which is 
sufficient to prove Theorem \ref{j+} and several other results from polynomial 
dynamics needed below. In the proof, we follow the ideas of \cite{p1}, 
with appropriate simplifications arising from this special case. 
For simplicity, we assume that the compact sets involved are infinite. 


We recall that, for a given compact set $K \subset \C$, a monic polynomial  
$P(z) \in \C[z]$ of degree $n > 0$ is called the $n$th polynomial of least 
deviation from zero if  
\[
\|P\|_K \leq \|Q\|_K
\]  
for any monic polynomial $Q(z) \in \C[z]$ of degree $n$, where  
\[
\|P\|_K := \max_{z \in K} |P(z)|.
\]  
It is well known that such a polynomial is unique whenever
$\operatorname{card} K \geq n$
(see, e.g., \cite[Chapter~1, Theorem~2.6]{ds}).
We denote it by $P_{n,K}$. 

The following lemma was proved in \cite{kamo} using the Kolmogorov
criterion for polynomials of least deviation (\cite{kol}; see also \cite[Chapter~1, Theorem~4.1]{ds}). For $n=1$, it
was also proved independently by a different method in the context of
equation \eqref{eqv} in \cite{p2}. Below we provide a proof following
the method of \cite{p1} (cf. Theorem 2.3 in \cite{p1}).

\bl \l{lemma}  
Let $K$ be an infinite compact set in $\C$ and let $X$ be a monic polynomial 
of degree $d\geq 1$. Then for any $n\geq 1$ the equality  
\[
 P_{n,K} \circ X=P_{nd,X^{-1}(K)} 
\]  
holds.  
\el  

\pr 
For any polynomial $Q$, define its averaging by $X$ by the formula 
\[
Q_X(z) = \frac{1}{d} \sum_{\substack{\zeta \in \C, \\ X(\zeta) = X(z)}} Q(\zeta),
\]
where each root $\zeta$ of multiplicity $k$ of \be \l{roo} X(\zeta) - X(z) = 0\ee is counted $k$ times. Clearly, we have
\be \l{1a}
\max_{z \in X^{-1}(K)} \Big| Q_X(z) \Big| 
\le \max_{z \in X^{-1}(K)} \frac{1}{d} \sum_{\substack{\zeta \in \C, \\ X(\zeta) = X(z)}} \Big| Q(\zeta) \Big| 
\le \max_{z \in X^{-1}(K)} \Big| Q(z) \Big|.
\ee

Observe now that for any polynomial $A$ of degree less than $d$, the function $A_X(z)$ is constant. The case $d = 1$ is trivial, so assume $d \ge 2$. For the monomials $A(z) = z^j,$  $1 \le j \le d-1$, the claim follows from the Newton formulas, which express $A_X(z)$ in terms of the elementary symmetric functions $S_j$, $1\leq j\leq d-1$,
of the roots of \eqref{roo}, which are independent of $z$
by Vieta's formulas. The general case then follows by linearity.  

The above implies that if $Q$ is a polynomial of degree $nd$ with $X$-adic decomposition
\be \l{in} 
Q(z) = \sum_{i=0}^{n} A_i(z) X^i(z),
\ee
then 
\be \l{bby} 
Q_X(z) = \sum_{i=0}^{n} a_i X^i(z),
\ee 
where $a_i \in \C$ for $0 \le i \le n$. Moreover, since $\deg Q = nd$, the degree of $A_n(z)$ in \eqref{in} is zero. Consequently, if $Q$ is monic, then $A_n(z)=a_n = 1$, and thus $Q_X(z)$ is monic.

Equality \eqref{bby} implies that, for any monic polynomial $Q$ of degree $dn$,  
\begin{equation} \label{2a}
\begin{split}
\max_{z \in X^{-1}(K)} \Big| Q_X(z) \Big| 
&= \max_{z \in K} \Big| \sum_{i=0}^{n} a_i z^i \Big| 
\ge \max_{z \in K} \Big| P_{n,K}(z) \Big| =\\
&= \max_{z \in X^{-1}(K)} \Big| (P_{n,K}\circ X)(z) \Big|.
\end{split}
\end{equation}
It then follows from \eqref{1a} and \eqref{2a} that
\[
\max_{z \in X^{-1}(K)} \Big| Q(z) \Big| 
\ge \max_{z \in X^{-1}(K)} \Big| (P_{n,K}\circ X)(z) \Big|.
\]
Therefore, the polynomial $P_{n,K} \circ X$ is the $nd$th polynomial of least deviation from zero on $X^{-1}(K)$. \qed

\bt \l{be} 
Let $A_1$, $A_2$ be non-constant polynomials and $K_1$, $K_2$ infinite compact sets in $\C$ satisfying  
\[
A_1^{-1}(K_1) \,=\, A_2^{-1}(K_2).
\] 
Assume that $\deg A_1 \mid \deg A_2$. 
Then there exists a polynomial $A$ such that 
\[
A_2 = A \circ A_1 \quad \text{and} \quad K_1 = A^{-1}(K_2). 
\]
\et 

\pr 
Let $a_1$ and $a_2$ be the leading coefficients of the polynomials $A_1$ and $A_2$, $n_1$ and $n_2$ their degrees, and let 
\[
K = A_1^{-1}(K_1) \,=\, A_2^{-1}(K_2). 
\]
Further, let $D$ be the disc of minimal radius containing $K_2$, and let $z_0 \in \C$ be 
its center. 

Replacing $K_2$ and $A_2$ by $K_2 - z_0$ and $A_2 - z_0$, without loss of generality we may 
assume that $z_0 = 0$, which implies $P_{1,K_2} = z$. Moreover, in this case, 
for the set 
$
\widehat{K_2} = (z/a_2)(K_2),
$ 
we also have $P_{1,\widehat{K_2}} = z$. Since  
\[
\bigl((z/a_2) \circ A_2\bigr)^{-1}(\widehat{K_2}) 
= A_2^{-1}(K_2) = K,
\] 
Lemma~\ref{lemma} applied with $n=1$ shows that the polynomial  
$
(z/a_2) \circ A_2
$ 
is the $n_2$-th polynomial of least deviation from zero on $K$.

On the other hand, setting 
$
\widehat{K_1} = (z/a_1)(K_1),
$ 
we have 
\[
K = (z / a_1 \circ A_1)^{-1}(\widehat{K_1}).
\] 
Thus, by Lemma \ref{lemma}, the polynomial 
\[
P_{n_2/n_1,\widehat{K_1}} \circ z / a_1 \circ A_1
\] 
is also the $n_2$-th polynomial of least deviation from zero on $K$.  
By the uniqueness of the polynomial of least deviation, we conclude that  
\[
P_{n_2/n_1,\widehat{K_1}} \circ z / a_1 \circ A_1 = z / a_2 \circ A_2,
\] 
which implies that the equality $A_2 = A \circ A_1$ holds for the polynomial 
\[
A = a_2 z \circ P_{n_2/n_1,\widehat{K_1}} \circ z / a_1.
\]

Finally, since \( A_2 = A \circ A_1 \), we have
\[
K =A_1^{-1}(K_1) \quad \text{and} \quad K= A_2^{-1}(K_2) = A_1^{-1}(A^{-1}(K_2)),
\]
which implies that
\[
A_1(K) = K_1 \quad \text{and} \quad A_1(K) = A^{-1}(K_2).
\]
Thus, \( K_1 = A^{-1}(K_2) \). \(\qed\)

\begin{example}\l{ex2}
Theorem~\ref{be} does not hold for arbitrary rational functions. This can be seen by considering the rational functions
\[
A_1 = \frac{z^2 - 1}{z^2 + 1} \quad \text{and} \quad A_2 = z + \frac{1}{z}
\]
of the same degree two, 
which appear in Example~\ref{ex1}. Indeed, it follows from \eqref{rr1} that \eqref{eqv} holds for 
\[
K_1 = \left(\frac{z^2 - 1}{z^2 + 1}\right)^{\!-1}\!\!(K), \qquad K_2 = \left(-\frac{2}{z^2 - 2}\right)^{\!-1}\!\!(K),
\]
where $K$ is any subset of $\C$, while equality \eqref{rr2} is impossible. 
\end{example}

\bt \l{j+} 
Let $A$ and $B$ be polynomials of degree at least two, and let $X$ be a non-constant polynomial. Then the relation
$A \underset{X}{\leq} B$ implies that 
\be \l{osn}
X^{-1}(J(A)) = J(B).
\ee 
Conversely, if for given $B$ and $X$ the condition \be \l{nos} X^{-1}(K) = J(B) \ee holds for some compact set $K \subset \C$, then there exists a polynomial $A$ such that \[ A \underset{X}{\leq} B \quad \text{and} \quad J(A) = K. \] \et

\pr 
The first part is known to hold for arbitrary rational functions $A$, $B$, and $X$ (see e.g. \cite{bu}, Lemma 5). 

Assume now that \eqref{nos} holds. Then  
\be
X^{-1}(K) = J(B) = (X \circ B)^{-1}(K).
\ee  
Applying now Theorem~\ref{be} with $A_1 = X$ and $A_2 = X \circ B$, we conclude that there exists a polynomial $A$ such that $A \underset{X}{\leq} B$.  
By the first part of the theorem, this  implies that \eqref{osn} holds.  
 Finally,  \eqref{osn} combined with \eqref{nos}   gives
$
J(A) = K.$ \qed

\begin{example}\l{ex3}
Theorem~\ref{j+} also fails for arbitrary rational functions. 
To see this, set 
\[
B = \frac{z^2 - 1}{z^2 + 1} \quad \text{and} \quad X = z + \frac{1}{z}.
\]
Then \eqref{rr1} implies that \eqref{nos} holds for 
\[
K = \left(-\frac{2}{z^2 - 2}\right)^{-1} \bigl(J(B)\bigr).
\]

Since $0$ belongs to the superattracting cycle $0\mapsto -1\mapsto 0$
of $B$, we have $0\notin J(B)$.
Therefore, since the rational map $z\mapsto -2/(z^2-2)$ sends $\infty$ to $0$,
the inverse image $K$ of $J(B)$ is compact in $\mathbb{C}$.
However, the equality \eqref{1} cannot hold for any rational function $A$. 
Indeed, let $F_1$ (resp. $F_2$) be the degree‑four rational function defined by either side of \eqref{rr1} (resp. \eqref{1}). Clearly, \eqref{rr1} and \eqref{1} imply that the two functions $F_1$ and $F_2$ belong to the field 
$$k = \mathbb{C}(B) \cap \mathbb{C}(X).$$
Thus, by L\"uroth's theorem, $k$ is generated by some non-constant rational function $F$, whose degree is either $2$ or $4$. Moreover, $\deg F \neq 2$. Indeed, otherwise both $X$ and $B$ would lie in $k$, so equality
\eqref{rr2} would hold, which is impossible. Therefore, each of the functions $F_1$ and $F_2$ generates $k$, implying that there exists a rational function $\mu$ of degree one such that $F_1 = \mu \circ F_2$. Thus,
\[
B \circ B = \mu \circ X \circ B,
\]
which again forces \eqref{rr2}.
\end{example}

\subsection{Polynomials with non-trivial $\Sigma(A)$ and polynomials sharing Julia sets} Let
\[
A(z) = a_0 z^d + a_1 z^{d-1} + \dots + a_d
\]
be a polynomial of degree $d\geq 2$. We say that $A$ is \emph{centered} if $a_1 = 0$. 
We recall that for a polynomial $A$, we defined $\Sigma(A)$ as the group of 
symmetries of $J(A)$, that is, the group of affine transformations 
$\mu(z) = az + b$ satisfying $\mu(J(A)) = J(A)$.

In this section, we recall two theorems concerning the Julia sets of polynomials.  
The first connects the nontriviality of $\Sigma(A)$ with the form of $A$,  
and the second provides a classification of polynomials that share a Julia set.  
These theorems follow from the results in \cite{b}, \cite{be2}, \cite{be1},  and \cite{sh}, but they can also be obtained using Theorems \ref{be} and \ref{j+}. Since these alternative proofs may be of independent interest, we present them here. 
Throughout, $\v_n$ denotes a primitive $n$th root of unity.

\bt \l{cic}  Let $A$ be a polynomial of degree $d \geq 2$ not conjugate to $z^d$.  
Then $\Sigma(A)$ is a finite cyclic group. Moreover, if $A$ is centered,  
then $\Sigma(A)$ is generated by the transformation $z \mapsto \v_n z$,  
where $n$ is the largest integer such that $A$ can be written in the form  
$
A = z^s R(z^n)
$
for some polynomial $R$ and  $s \geq 0$.  

\et  

\pr Without loss of generality, we may assume that $A$ is centered. If $\mu \in \Sigma(A)$, then  
\[
A^{-1}(J(A)) = (A \circ \mu)^{-1}(J(A)).
\]  
Hence, by Theorem~\ref{be}, the equality
\be \l{no}
A \circ \mu = \nu \circ A
\ee
holds for some degree-one polynomial $\nu$, which also belongs to $\Sigma(A)$. 

Since $A$ is centered, \eqref{no} implies that $\mu(0)=0$. Moreover, since this equality holds for every $\mu\in\Sigma(A)$, taking $\mu=\nu$ gives $\nu(0)=0$.
Now \eqref{no} implies that, unless $A$ is conjugate to a power,
$\mu=\v z$ for some root of unity $\v$. Furthermore, if $\v$ has order $n$, then $A$ is of the form $A = z^s R(z^n)$.

Conversely, if $A$ has the form $z^s R(z^n)$, then considering the semiconjugacy diagram 
\[
\begin{CD}  
\mathbb{CP}^1 @> z^s R(z^n) >> \mathbb{CP}^1 \\  
@V z^n VV @VV z^n V \\  
\mathbb{CP}^1 @> z^s R^n(z) >> \mathbb{CP}^1  
\end{CD}  
\]  
and applying Theorem~\ref{j+}, we obtain the equality  
\[
J(A) = J(z^s R(z^n)) = (z^n)^{-1} \bigl( J(z^s R^n(z)) \bigr),
\]  
which shows that $J(A)$ is invariant under $z \mapsto \v_n z$.  \qed

Theorem \ref{cic} implies the following well-known fact, which we will often use without mentioning it explicitly.  

\bc Let $A$ be a polynomial of degree at least two, and let $\mu \in \Sigma(A)$.  
Then  
 $J(\mu \circ A) = J(A)$.
\ec 
\pr 
Without loss of generality, we may assume that $A$ is centered and not conjugate to a power, since in the latter case the corollary is obviously true.
 
Let $n = |\Sigma(A)|$. Then, by Theorem~\ref{cic}, 
$A = z^s R(z^n)$ and $\mu = \v z$, where $\v$ is an $n$th root of unity. Thus, $z^n \circ \mu = z^n$, implying that the relations  
\[  
z^s R^n(z) \underset{z^n}{\leq} A \quad \text{and} \quad z^s R^n(z) \underset{z^n}{\leq} \mu \circ A  
\]  
hold simultaneously.  
Hence, $$J(\mu \circ A) = J(A)=(z^n)^{-1}J(z^s R^n(z))$$ by the first part of Theorem~\ref{j+}. 
 \qed  

\bc \l{lup}
Let $A$ be a polynomial of degree $d\geq2$ not conjugate
to a power, and let $\mu\in\Sigma(A)$ have order $m$.
Then $A$ and $\mu\circ A$ share an iterate if and only if
$\gcd(d,m)=1$.
\ec
\pr
Let $n=|\Sigma(A)|$.
As above, we may assume that
$
A = z^s R(z^n)
$
and $\mu = \v z$, where $\v$ is an $n$th root of unity of order $m$.
Since $d\equiv s\pmod n$,
for any integer $N\geq1$ we have
\be \l{twistit}
(\mu \circ A)^{\circ N} = \v^{\,1+d+d^2+\dots+d^{\,N-1}} \, A^{\circ N},
\ee
implying that if $A$ and $\mu\circ A$ share an iterate, then necessarily $\gcd(d,m)=1$.

Conversely, if $\gcd(d,m)=1$, then the map
$\kappa:t\mapsto dt+1$ is a permutation of $\mathbb Z/m\mathbb Z$.
Let $N$ be the length of the cycle containing $0$. Then
\[
1+d+\cdots+d^{N-1}\equiv \kappa^{\circ N}(0)\equiv0\pmod m.
\]
By \eqref{twistit}, this implies
$$(\mu\circ A)^{\circ N}=A^{\circ N}.\eqno{\Box}$$

The following result provides a classification of polynomials sharing a Julia set.

\bt \l{sha}
Let $K \subset \C$ be a compact set that is the Julia set of a non-special polynomial.
Then there exists a polynomial $Q$ such that $J(Q) = K$, and any polynomial $P$ with $J(P) = K$ has the form
$
P = \mu \circ Q^{\circ n}
$
for some $\mu \in \Sigma(Q)$ and $n \geq 1$.
\et

\pr  
Let $Q$ be a non-special polynomial of minimal degree with $J(Q) = K$. 
We will show that for $Q$ defined in this way, the conclusion of the theorem 
holds for any polynomial $P$  of degree at least two with \be \l{bbs} P^{-1}(K) = K.\ee 

We begin by showing that \eqref{bbs} implies that $\deg Q$ divides $\deg P$. 
Since  
\[
(P \circ Q^{\circ i})^{-1}(K) = (Q^{\circ i} \circ P)^{-1}(K) = K, \quad i\geq 1,
\] 
 it follows from Theorem~\ref{be} that for every $i \geq 1$ the equality
\be \l{eqref} 
Q^{\circ i} \circ P = \mu_i \circ P \circ Q^{\circ i}
\ee  holds
for some $\mu_i \in \Sigma(Q)$.  
Furthermore, since $Q$ is not conjugate to a power, the group $\Sigma(Q)$ is finite, implying that there exist $i_2 > i_1$ such that $\mu_{i_1} = \mu_{i_2}$. Setting now $\mu = \mu_{i_1} = \mu_{i_2}$, we see that \eqref{eqref} yields 
\[
Q^{\circ (i_2 - i_1)} \circ \mu \circ P \circ Q^{\circ i_1} = \mu \circ P \circ Q^{\circ i_2},
\]
which implies that $\mu \circ P$ commutes with $Q^{\circ (i_2 - i_1)}$.

Recall that by Ritt’s theorem (see \cite{r}, and also \cite{e2}, \cite{rev}), if rational functions $A$ and $B$ of degree at least two commute, then either they share an iterate, or they are both Latt\`es maps, or they are both conjugate to powers (with positive or negative exponents), or they are both conjugate to Chebyshev polynomials (with positive or negative signs).  
Since a polynomial cannot be a Latt\`es map, and $ Q$ is non-special, we conclude that 
\be \l{ded}  
 Q^{\circ l} = (\mu \circ P)^{\circ k}  
\ee  
for some  $k, l \geq 1$, which implies  
\[
(\deg Q)^l = (\deg P)^k.  
\]  
Since by assumption $\deg Q \leq \deg P$, we have $l \geq k$, which yields  that  $\deg Q$ divides  $\deg P$.  

Let now $P$ be an arbitrary polynomial of degree at least two  such that \eqref{bbs} holds.  
Then it follows from  
\[
P^{-1}(K) = Q^{-1}(K) = K
\]  
by Theorem~\ref{be} that  
$  
P = \t P \circ Q  
$  
for some polynomial $\t P$ with   
$
\t P^{-1}(K) = K.  
$  
Repeating this argument as long as the remaining factor
has degree at least two, we eventually obtain the representation
$$
P = \mu \circ Q^{\circ n},  
$$  
where  $\mu\in \Sigma(Q)$  and $n \geq 1$. 
\qed  

\bc \l{shaa} 
Let $P$ and $Q$ be polynomials of the same degree $d \geq 2$ such that $J(P) = J(Q)$.  
Then $Q = \mu \circ P$ for some $\mu \in \Sigma(P)$.  
\ec  

\pr
For non-special polynomials $P$ and $Q$, this follows from Theorem~\ref{sha}, and it is straightforward to verify that the same conclusion holds when one or both of them are special. Alternatively, the corollary is an immediate consequence of Theorem~\ref{be}.
\qed

\section{\l{s4} Semiconjugacies involving iterates}
In this section, we analyze the consequences of the relation $A \underset{X}{\leq} B$ for polynomials $A$ and $B$ in the case where one of them is an iterate of some polynomial. 

We begin with the following result.

\bt \l{t30}  
Let $A$ and $B$ be polynomials of degree at least two, and $X$ a non-constant polynomial such 
that  $A \underset{X}{\leq} B$.  Then, for every polynomial $\widehat{B}$ with \linebreak $J(\widehat{B}) = J(B)$, there exists a polynomial $\widehat{A}$ such that  
$\widehat{A} \underset{X}{\leq} \widehat{B}$ and $J(\widehat{A}) = J(A)$.

\et  

\pr  
By the first part of Theorem~\ref{j+},  
$$
X^{-1}\!(J(A)) = J(B).
$$  
Since $J(\h B) = J(B)$, this implies, by the second part of Theorem~\ref{j+}, that there exists a polynomial  
$\widehat{A}$ such that  
$\widehat{A} \underset{X}{\leq} \h B$ and $J(\widehat{A}) = J(A)$.  
 \qed  

Theorem \ref{t30} implies the following result. 

\bt \l{mt1}  
Let $A$ and $B$ be polynomials of degree at least two, and $X$ a non-constant polynomial such 
that  $A \underset{X}{\leq} B$.  Suppose that  
$B = \hat{B}^{\circ l}$ for some polynomial $\hat{B}$ and some  $l \geq 2$.  
Then there exists a polynomial $\h{A}$ such that $\h A   \underset{X}{\leq} \h B$ and $A = \h{A}^{\circ l}$. \et  

\pr  Since $J (\hat{B})=J(B)$,  Theorem \ref{t30} implies that   
there exists a polynomial $\h{A}$ such that $\h A   \underset{X}{\leq} \h B$. In turn, this relation gives 
$
\widehat{A}^{\circ l} \underset{X}{\leq} \widehat{B}^{\circ l},
$ 
which, together with
$
A \underset{X}{\leq} \widehat{B}^{\circ l},
$ 
implies that 
$
A = \widehat{A}^{\circ l}.
$ 
\qed

\bc\l{43} 
Let $A$ and $B$ be polynomials of degree at least two such that $A \sim B$. 
Then $A$ is an $l$th iterate of some polynomial if and only if $B$ is an $l$th iterate of some polynomial.
\ec

\pr
Since the relation $A \sim B$ implies the relations $A \leq B$ and $B\leq A$, the statement  follows  from Theorem~\ref{mt1}. \qed

\bl \l{spec}  
Let \( A \) be a polynomial of degree \( d > 2 \). Assume that the equality  
\( A^{\circ l} = \mu \circ T_{d^l} \circ \nu \)  
holds for some \( l \geq 2 \) and some polynomials \( \mu \) and \( \nu \) of degree one.  
Then  $A^{\circ l}$ is conjugate to $\pm   T_{d^l}.$ 
\el  

\pr  
Clearly, without loss of generality, we may assume that \( \nu = \operatorname{id} \), that is,  
\be \l{am} 
A^{\circ l} = \mu \circ T_{d^l}, 
\ee  
which, by Theorem~\ref{r1}, implies that
\be \l{by}
A = \mu \circ T_d \circ \alpha \quad {\rm and} \quad A^{\circ (l - 1)} = \alpha^{-1} \circ T_{d^{l - 1}}
\ee  
for some degree one polynomial \( \alpha \). 

It follows from \eqref{am} and \eqref{by} that  
\[
\mu \circ T_{d^l} = A^{\circ (l - 1)} \circ A = \alpha^{-1} \circ T_{d^{l - 1}} \circ \mu \circ T_d \circ \alpha,
\]  
so  
\[
T_{d^l} = (\mu^{-1} \circ \alpha^{-1} \circ T_{d^{l - 1}}) \circ (\mu \circ T_d \circ \alpha),
\]  
which, again by Theorem~\ref{r1}, implies that 
\be \l{sid}
\mu^{-1} \circ \alpha^{-1} \circ T_{d^{l - 1}} = T_{d^{l - 1}} \circ \beta \quad {\rm and} \quad
 \mu \circ T_d \circ \alpha=\beta^{-1} \circ T_d 
\ee  
for some degree one polynomial \( \beta \).

Since, as shown in the proof of Theorem~\ref{i}, the equality 
\be \l{pos} 
T_d = \gamma \circ T_d \circ \delta, \quad d > 2,
\ee  
for some degree one polynomials \( \gamma \) and \( \delta \) implies that  
$\delta = \pm z$ and $\gamma = \pm z,$ 
we obtain  
\be \alpha \circ \mu = \pm z, \quad 
\beta = \pm z, \quad \beta \circ \mu = \pm z, \quad
\alpha = \pm z.
\ee  
Hence $\mu = \pm z$, and $A^{\circ l} = \pm T_{d^l}$ by \eqref{am}. 
\qed

The following lemma is well-known. For the sake of completeness, we provide a sketch of the proof.

\bl \l{gop}  
Let $A$ be a polynomial of degree $d \geq 2$, and let $n > 1$ be an integer. Then there exist at most two complex numbers $\alpha$ such that all but one of the multiplicities of $A$ at the points of $A^{-1}(\alpha)$ are divisible by $n$. Furthermore, if there exist two such numbers, then $n = 2$ and there exist degree-one polynomials $\mu$ and $\nu$ such that  
$
A = \mu \circ T_d \circ \nu.
$
\el  

\pr  
It follows from the Riemann--Hurwitz formula that the preimage of a set $K \subset \C$ containing $k$ elements under $A$ contains at least $(k-1)d+1$ points, with equality if and only if $K$ contains all finite critical values of $A$. Applying this to the case where \(K=\{\alpha_1,\alpha_2,\dots,\alpha_k\}\) is the set of all complex numbers satisfying the condition of the lemma, we obtain  
\be \l{x} 
(k-1)d + 1 \leq k\left(\frac{d-1}{n} + 1\right),
\ee  
which implies  
\[
k(d-1)\left(1 - \frac{1}{n}\right) \leq d-1 \quad \text{and} \quad k\left(1 - \frac{1}{n}\right) \leq 1.
\]  
Thus,
$
k \le 2,
$
and if \(k = 2\), then \(n = 2\) and the inequality in \eqref{x} becomes an equality.
 This implies that
$\alpha_1$ and $\alpha_2$ are the only finite critical values of $A$, and that all
multiplicities of $A$ at the points in $A^{-1}\{\alpha_1,\alpha_2\}$ are two,
except at two points whose multiplicities are one.

Finally, it is well known that in the last case 
$
A = \mu \circ T_d \circ \nu
$
for some degree-one polynomials $\mu$ and $\nu$
(see, e.g., \cite[Lemma~4, p.~27]{sch}). 
An elegant way to see this is via the correspondence between equivalence classes of polynomial Belyi functions and plane trees, which is a special case of Grothendieck’s theory of dessins d’enfants (see, e.g., \cite{lz}). Under this correspondence, the statement reduces to the fact that the only trees whose vertex valencies do not exceed two are ``chains.'' \qed  

\smallskip

 The following result is an analogue of Theorem~\ref{mt1}
with the roles of $A$ and $B$ reversed.
Unlike Theorem~\ref{mt1}, it gives a conclusion for an iterate
of $B$ rather than for $B$ itself.

\bt \label{t33}
Let $A$ and $B$ be polynomials of degree at least two such that
$A\leq B$. Suppose that $A={\h A}^{\circ l}$ for some
polynomial $\h A$ and some $l\geq2$.
Then there exists a polynomial $\widehat B$ such that
$B^{\circ e}={\h B}^{\circ le}$ for some $e\geq1$.
\et
\pr  
If $A$ is special, the claim follows from Theorem \ref{spe}, taking into account that for even $d$ the polynomials $T_d$ and $-T_d$ are conjugate,  while for $d$ odd their second iterates are equal. Thus, it suffices to consider the case where $A$, and hence $B$, are  non-special.

Let $X$ be a polynomial such that $A \underset{X}{\leq} B$. By Lemma~\ref{25}, we can find polynomials  $ B_0$, $U$, and $X_0$ such that $B\sim B_0$, $X=X_0\circ U$, the diagram  
\[ \l{retu} 
\begin{CD}
\mathbb{C}\mathbb{P}^1 @> B >> \mathbb{C}\mathbb{P}^1 \\
@V U VV @VV U V \\
\mathbb{C}\mathbb{P}^1 @> {B}_0 >> \mathbb{C}\mathbb{P}^1 \\
@V X_0 VV @VV X_0 V \\
\mathbb{C}\mathbb{P}^1 @>  {\h A}^{\circ l} >> \mathbb{C}\mathbb{P}^1
\end{CD}
\]
commutes, and 
 \(  {\h A}^{\circ l}, B_0, X_0  \) is a primitive solution of equation~\eqref{1}. Furthermore, we may assume that $\deg X_0 > 1$, since otherwise ${\hat A}^{\circ l} \sim B$, and the statement of the theorem follows from Corollary~\ref{43}.  
Since $A$ is not special, it follows from Theorem~\ref{i} that, without loss of generality, we may assume that
\[
X_0 = z^n, \quad B_0 = z^s R(z^n), \quad {\h A}^{\circ l} = z^s R^n(z),
\]  
for some non-constant polynomial $R$, with $n \geq 2$, $s \geq 1$, and $\gcd(s, n) = 1$.

We may also assume that $\deg\widehat A>2$. 
Indeed, if $\deg\widehat A=2$, then $\deg A=2^l$,
so $n$ is odd, since $\gcd(n,\deg A)=1$.
As $A$ is non-special, $R$ has a nonzero root.
The multiplicity of $A=z^sR^n(z)$ at this root is divisible
by $n>1$, whereas every multiplicity of an iterate
of a quadratic polynomial is a power of two.
This is a contradiction.

To prove the theorem, it is enough to show that there exists a polynomial 
$\widehat{B}_0$ such that the equality 
\be \l{tak} 
 B_0^{\circ e} =\widehat{B}_0^{\circ le}  
\ee  
holds for some integer $e\geq1$. Indeed, since $B \sim B_0$ implies 
$B^{\circ e} \sim B_0^{\circ e}$,  
Corollary~\ref{43} then implies that 
$$B^{\circ e} = \widehat{B}^{\circ le}$$ for some polynomial $\widehat{B}$.

Applying Theorem~\ref{r1} to the equality
\[
{\h A} \circ ({\h A}^{\circ (l - 1)} \circ X_0) = X_0 \circ B_0,
\]  
we can find polynomials  $B_1$,  $B_2$, and $Y$ such that $B_0 = B_1 \circ B_2$, and the diagram  
\be \label{mi}  
\begin{CD}
\mathbb{CP}^1 @> B_2 >> \mathbb{CP}^1 @> B_1 >> \mathbb{CP}^1 \\
@V X_0 VV @V Y VV @V X_0 VV \\
\mathbb{CP}^1 @> {\h A}^{\circ (l - 1)} >> \mathbb{CP}^1 @> {\h A} >> \mathbb{CP}^1
\end{CD}
\ee  
commutes.   
Moreover, since $n = \deg X_0$ is coprime to $$\deg B_0 = \deg B=\deg A,$$ setting $d = \deg {\h A}$, we have
\begin{multline}
\deg B_2 = \gcd\bigl(\deg B_0,\; \deg({\h A}^{\circ (l-1)} \circ X_0)\bigr)= \\
= \gcd(d^l,\; d^{l-1} \deg X_0) = d^{l-1} = \deg {\h A}^{\circ (l-1)}.
\end{multline}
Thus, 
\[
\deg Y = \deg X_0 = n \quad \text{and} \quad \deg B_1 = \deg {\h A}.
\]

The commutativity of \eqref{mi} implies that the diagram
\be 
\begin{CD}
\mathbb{CP}^1 @> B_1 >> \mathbb{CP}^1 @> B_2 >> \mathbb{CP}^1 \\
@V Y VV @V X_0 VV @V Y VV \\
\mathbb{CP}^1 @> {\h A} >> \mathbb{CP}^1 @> {\h A}^{\circ (l - 1)} >> \mathbb{CP}^1
\end{CD}
\ee  
also commutes, which in turn yields the commutativity of  
\[
\begin{CD}
\mathbb{CP}^1 @> B_2 \circ B_1 >> \mathbb{CP}^1 \\
@V Y VV @VV Y V \\
\mathbb{CP}^1 @> {\h A}^{\circ l} >> \mathbb{CP}^1. 
\end{CD}
\]  
Since $\deg Y = n$, applying Theorem~\ref{i} again we conclude that there exist degree-one polynomials $\mu$ and $\nu$ such that  
\[
{\h A}^{\circ l} = \nu \circ z^{\tilde{s}} \tilde{R}^n(z) \circ \nu^{-1}, \quad 
Y = \nu \circ z^n \circ \mu, \quad 
B_2 \circ B_1 = \mu^{-1} \circ z^{\tilde{s}} \tilde{R}(z^n) \circ \mu,
\]  
for some polynomial $\tilde{R}$, with $\tilde{s} \geq 1$ and $\gcd(\tilde{s}, n) = 1$.   

Furthermore, we have  
\be \label{v} 
\nu(z) = a z, \quad a \in \mathbb{C}^*.
\ee  
Indeed, if $\nu(0)  \ne 0$, then from  
\[
{\h A}^{\circ l} = z^s R^n(z) = \nu \circ z^{\tilde{s}} \tilde{R}^n(z) \circ \nu^{-1},
\]  
it follows that the conditions of Lemma~\ref{gop} are satisfied for $\alpha= \nu(0)$ and $\alpha = 0$, implying  that 
\[
{\h A}^{\circ l} = \mu \circ T_m \circ \nu
\]  
for some degree-one polynomials $\mu$ and $\nu$. Hence, $A={\h A}^{\circ l} $
 is special by Lemma~\ref{spec},  in contradiction with the assumption.

Since $X_0=z^n$, it follows from \eqref{v} that 
\[
Y = \nu \circ z^n \circ \mu = z^n \circ (\widehat{a} z) \circ \mu,
\]  
where $\widehat{a}^n = a$. Thus, the commutativity of the right square in \eqref{mi} implies the commutativity of the diagram  
\be \label{b0} 
\begin{CD}
\mathbb{CP}^1 @> \widehat{B}_0 >> \mathbb{CP}^1 \\
@V z^n VV @VV z^n V \\
\mathbb{CP}^1 @> {\h A} >> \mathbb{CP}^1
\end{CD}
\ee  
where  
\[
\widehat{B}_0 = B_1  \circ \mu^{-1}\circ (\widehat{a}^{-1} z).
\]  
Thus, 
$
{\h A} \underset{z^n}{\leq} \widehat{B}_0$ and hence \be \l{rrt} {\h A}^{\circ l} \underset{z^n}{\leq} \widehat{B}_0^{\circ l}.\ee  

Combined with  
$
{\h A}^{\circ l} \underset{z^n}{\leq} B_0,
$  
the relation \eqref{rrt} 
 yields  
\[
z^n \circ \widehat{B}_0^{\circ l} = z^n \circ B_0,
\]  
implying that  
\be \label{from} 
\widehat{B}_0^{\circ l}= \mu_0 \circ B_0,
\ee  
where $\mu_0 = \v z$ with $\v^n = 1$. 
Finally,  $\mu_0\in\Sigma(B_0)$ by Theorem~\ref{cic}, and 
the order of $\mu_0$  is coprime
to $\deg B_0$, since $$\deg B_0\equiv s\pmod n$$ and
$\gcd(s,n)=1$. 
By Corollary~\ref{lup} and \eqref{from}, we obtain
\[
B_0^{\circ e}
=(\mu_0\circ B_0)^{\circ e}
=\widehat B_0^{\circ le}
\]
for some integer $e\geq1$.
\qed

\section{Proof of Theorem \ref{t1}} 

The following result is an algebraic counterpart of Theorem~\ref{t1}.

\bt \l{fin0}
Let $A_1$ and $A_2$ be polynomials of degree $d\geq2$
not conjugate to a power, and $X_1$ and $X_2$ non-constant
polynomials satisfying
\[
A_1^{\circ k}\underset{X_1}{\leq}B,
\quad
A_2^{\circ k}\underset{X_2}{\leq}B
\]
for some polynomial $B$ and $k\geq1$.
Then there exist a polynomial $\widehat B$ and
$\gamma_1\in\Sigma(A_1)$, $\gamma_2\in\Sigma(A_2)$ such that
\be \l{ree1}
\gamma_1\circ A_1\underset{X_1}{\leq}\widehat B,
\quad
\gamma_2\circ A_2\underset{X_2}{\leq}\widehat B.
\ee
Moreover, there exists an integer
\[
1\leq L\leq\operatorname{lcm}\bigl(|\Sigma(A_1)|,|\Sigma(A_2)|\bigr)
\]
such that
\be \l{ree2}
A_1^{\circ L}\underset{X_1}{\leq}\widehat B^{\circ L},
\quad
A_2^{\circ L}\underset{X_2}{\leq}\widehat B^{\circ L}.
\ee
\et
\pr
Applying Theorem~\ref{t33} to either of the relations
\be \l{re0}
A_1^{\circ k}\underset{X_1}{\leq}B,
\quad
A_2^{\circ k}\underset{X_2}{\leq}B,
\ee
we obtain an integer $e\geq1$ and a polynomial $\widehat B$
such that $B^{\circ e}=\widehat B^{\circ ek}$, implying that 
\be \l{is}
A_1^{\circ ek}\underset{X_1}{\leq}\widehat B^{\circ ek},
\quad
A_2^{\circ ek}\underset{X_2}{\leq}\widehat B^{\circ ek}.
\ee
Therefore, by Theorem~\ref{t30}, there exist polynomials
$\widehat A_1$ and $\widehat A_2$ such that
\be \l{iss}
\widehat A_1\underset{X_1}{\leq}\widehat B,
\quad
\widehat A_2\underset{X_2}{\leq}\widehat B,
\ee
and
$$J(\widehat A_1)=J(A_1),\quad J(\widehat A_2)=J(A_2).$$
Moreover, since \eqref{is} and \eqref{iss} imply that
$$\deg A_1=\deg\widehat A_1,\quad\deg A_2=\deg\widehat A_2,$$
it follows from Corollary~\ref{shaa} that there exist
$\gamma_1\in\Sigma(A_1)$ and $\gamma_2\in\Sigma(A_2)$ such that
$$\widehat A_1=\gamma_1\circ A_1,\quad
\widehat A_2=\gamma_2\circ A_2.$$
Thus, the relations \eqref{ree1} hold.

Further, the relations
\be \l{si0}
\gamma_1\circ A_1\underset{X_1}{\leq}\widehat B,
\quad
\gamma_2\circ A_2\underset{X_2}{\leq}\widehat B
\ee
yield the relations
\be
(\gamma_1\circ A_1)^{\circ ek}\underset{X_1}{\leq}\widehat B^{\circ ek},
\quad
(\gamma_2\circ A_2)^{\circ ek}\underset{X_2}{\leq}\widehat B^{\circ ek}.
\ee
Together with \eqref{is}, this shows that
$\gamma_1\circ A_1$ shares an iterate with $A_1$, and
$\gamma_2\circ A_2$ shares an iterate with $A_2$.
Let $m_i$ be the order of $\gamma_i$, $i=1,2$.
Corollary~\ref{lup} gives $\gcd(d,m_i)=1$ for $i=1,2$.

Set $M=\operatorname{lcm}(m_1,m_2)$.
Since $\gcd(d,M)=1$, the map $\kappa:t\mapsto dt+1$ is a permutation
of $\mathbb Z/M\mathbb Z$, and the cycle containing $0$ has some
length $L\leq M$. Thus,
\[
1+d+\cdots+d^{L-1}\equiv \kappa^{\circ L}(0)\equiv0\pmod M.
\]
It follows from \eqref{twistit} that
\[
(\gamma_i\circ A_i)^{\circ L}=A_i^{\circ L},
\qquad i=1,2.
\]
Together with \eqref{si0}, this proves \eqref{ree2}.
Finally, since $m_i$ divides $|\Sigma(A_i)|$, we have
$$
L\leq M\leq\operatorname{lcm}\bigl(|\Sigma(A_1)|,|\Sigma(A_2)|\bigr). \eqno{\Box} 
$$

\vskip 0.2cm
\noindent{\it Proof of Theorem~\ref{t1}.}
We recall that if $A_1$ and $A_2$ are non-special polynomials
of degree at least two, then any irreducible
$(A_1,A_2)$-invariant curve $C$ that is neither a vertical
nor a horizontal line has genus zero and admits a generically
one-to-one parametrization by polynomials $X_1$ and $X_2$
such that the diagram
\be \l{dur}
\begin{CD}
(\C\P^1)^2 @>(B,B)>> (\C\P^1)^2 \\
@V (X_1,X_2) VV @VV (X_1,X_2) V\\
(\C\P^1)^2 @>(A_1,A_2)>> (\C\P^1)^2
\end{CD}
\ee
commutes for some polynomial $B$ (see Proposition~2.34
of \cite{ms} or Section~4.3 of \cite{pj}). Conversely,
if diagram~\eqref{dur} commutes, then the curve parametrized
by \linebreak $t\mapsto(X_1(t),X_2(t))$ is invariant under $(A_1,A_2)$.

If $C$ is periodic, applying this characterization to
an iterate of $(A_1,A_2)$ gives
\be \l{app}
A_1^{\circ k}\underset{X_1}{\leq}B,
\quad
A_2^{\circ k}\underset{X_2}{\leq}B
\ee
for some $k\geq1$ and polynomials $X_1,X_2,B$.
By Theorem~\ref{fin0}, there exist a polynomial $\widehat B$
and an integer
\[
1\leq L\leq\operatorname{lcm}\bigl(|\Sigma(A_1)|,|\Sigma(A_2)|\bigr)
\]
such that
\[
A_1^{\circ L}\underset{X_1}{\leq}\widehat B^{\circ L},
\quad
A_2^{\circ L}\underset{X_2}{\leq}\widehat B^{\circ L}.
\]
Hence $(A_1,A_2)^{\circ L}(C)=C$, and the period of $C$
is at most $L$.
\qed

\begin{example}\l{sharpperiod}
The bound in Theorem~\ref{t1} is attained for every $n\geq3$
by the polynomials
\[
A_1(z)=z^{n+1}+z,\qquad A_2(z)=\zeta A_1(z),
\]
where $\zeta$ is a primitive $n$th root of unity.
Indeed, each polynomial has $n$ distinct finite critical
values, so neither is special, and, by Theorem~\ref{cic},
\[
|\Sigma(A_1)|=|\Sigma(A_2)|=n.
\]
Let $C_j$ be the line $y=\zeta^j x$.  Since
$$A_1(\zeta^j z)=\zeta^j A_1(z),$$ we have
\[
(A_1,A_2)(C_j)=C_{j+1}.
\]
Thus the diagonal $C_0$ has period  $n$, attaining
the bound $\operatorname{lcm}(n,n)=n$.
\end{example}

Note that although the algebraic part of the proof of
Theorem~\ref{t1} only requires that $A_1$ and $A_2$ are not
conjugate to powers, the theorem does not extend to Chebyshev
polynomials. The reason is that while invariant curves are
still characterized by \eqref{dur}, the functions in this
diagram may be rational rather than polynomial.

\begin{example}
Fix $d\geq2$. For each integer $r\geq3$, let $\zeta$ be
a primitive root of unity of order $d^r-1$, and let $C_\zeta$
be the irreducible curve parametrized by
\[
x=\pi(t),\qquad y=\pi(\zeta t),
\qquad \pi(t)=\frac{t+t^{-1}}{2}.
\]
The identity $$T_d\circ\pi=\pi\circ z^d$$ gives
\[
(T_d,T_d)(C_\zeta)=C_{\zeta^d}.
\]
Since $C_\alpha=C_\beta$ if and only if
$\beta=\alpha^{\pm1}$, the period of $C_\zeta$ is the least
positive integer $k$ such that
\[
d^k\equiv\pm1\pmod{d^r-1},
\]
which is $r$. Thus, periods are unbounded for fixed
$(T_d,T_d)$.
\end{example}

\section{Semiconjugacies between iterates and proof of Theorem~\ref{t2}}
In this section we prove Theorems~\ref{fin00} and \ref{fin},
which allow us to control polynomial solutions of
$$
A^{\circ l} \underset{X}{\le} B^{\circ s}, \qquad s, l \ge 1,
$$
when one of the polynomials $A$ or $B$ is fixed, while the other,
together with the integers $s$ and $l$, may vary. Then we deduce Theorem~\ref{t2} 
 from these results.

We begin with the following two theorems.

\bt \label{ko1}  
For any $d \ge 2$, there exists a constant $r(d)$ such that 
for every polynomial $B$ of degree $d$, there are 
polynomials $A_1, A_2, \dots, A_r$ with $r \le r(d)$ 
satisfying the following property: for a polynomial $A$, the relation 
 $A \le B$ implies that $A$ is conjugate to one of 
$A_1, A_2, \dots, A_r$.
 
\et  

\pr If $B$ is special, the claim follows from Theorem~\ref{spe}. For non-special $B$, Theorem~\ref{ko1} is a particular case of Theorem~1.1 in \cite{fin}, which establishes the result in the broader context of rational functions. In the polynomial case, Theorem~\ref{ko1} also follows from Theorem~1.6 in \cite{pj}. 

Indeed, the latter theorem provides a polynomial $X_B$, whose
degree is bounded in terms of $d$, such that $A\le B$
implies that $A\underset{X}{\le}B$
for some compositional right factor $X$ of $X_B$.
By Theorem~\ref{r1}, right factors of $X_B$
of the same degree differ by postcomposition with a
degree-one polynomial and therefore yield conjugate
polynomials $A$. Therefore, the number of conjugacy classes
is at most $\deg X_B$.
\qed

\bt \l{k02} 
For any $d \ge 2$, there exists a constant $r(d)$ such that
for every polynomial $A$ of degree $d$, there are 
polynomials $B_1, B_2, \dots, B_r$ with $r \le r(d)$
satisfying the following property: for a polynomial $B$, the relation 
$A \le B$ implies that $B$ is conjugate to one of
$B_1, B_2, \dots, B_r$.
 
\et 

\pr If $A$ is special, the theorem follows from Theorem~\ref{spe}, so assume that $A$ is not special. 
It follows from Lemma~\ref{25} that to prove the theorem it suffices to establish the following two statements.  

First, for any $d \geq 2$, there exists a constant $r_1(d)$ such that for every polynomial $A$ of degree $d$, there exist polynomials $C_1, C_2, \dots, C_r$ with $r \leq r_1(d)$ such that if $A, B_0, X_0$ with $\deg X_0 \geq 2$ form a primitive solution of \eqref{1}, then $B_0$ is conjugate to one of $C_1, C_2, \dots, C_r$.  

Second, for any $d \geq 2$, there exists a constant $r_2(d)$ such that for every polynomial $C$ of degree $d$, there exist polynomials $B_1, B_2, \dots, B_r$ with $r \leq r_2(d)$ such that $C \sim B$ implies that $B$ is conjugate to one of $B_1, B_2, \dots, B_r$.  

The second statement follows from Theorem~\ref{ko1}, since $C \sim B$ implies $B \leq C$. Let us prove the first statement.
 If $A$ is non-special and $A, B_0, X_0$ with $\deg X_0 \geq 2$ form a primitive solution of \eqref{1}, then by Theorem~\ref{i} there exist degree-one polynomials $\mu$ and $\nu$ such that  
\be \l{fi} 
A = \nu \circ z^s R^n(z) \circ \nu^{-1}, \quad  
X_0 = \nu \circ z^n \circ \mu, \quad  
B_0 = \mu^{-1} \circ z^s R(z^n) \circ \mu,
\ee  
where $R$ is a non-constant polynomial, $n \geq 2$, $s \geq 1$, and $\gcd(s, n) = 1$. Furthermore, since $R$ is non-constant, we have $n < d$, so $n$ can take only finitely many values among different solutions. In addition, for each fixed $n$, Lemma~\ref{gop}, applied to the fixed
polynomial $A$, shows that $\nu(0)$ can take at most two
distinct values among these solutions.  
Hence, it suffices to prove the first statement for fixed $n \geq 2$ and fixed $\nu(0)$.

Consider a primitive solution of \eqref{1} of the form \eqref{fi},
and let $A,$ $\t X_0,$ $\t B_0$ be another primitive solution
satisfying the above conditions.
Then 
$$\t X_0 = \t\nu \circ z^{ n} \circ \t\mu$$ 
for some degree-one polynomials $\t\mu, \t\nu$. Furthermore, since $\t\nu(0)=\nu(0)$ by assumption, 
$$\t X_0 
= X_0 \circ \delta$$
for some degree-one polynomial $\delta$.
Thus, $$A\circ \t X_0= \t X_0\circ \t B_0$$
implies that 
\be \l{then} A \circ  X_0 = X_0\circ \delta \circ \t B_0 \circ \delta^{-1}.\ee

Let $\mu_i$, $1 \leq i \leq n$, be all M\"obius transformations satisfying  
$$X_0 \circ \mu = X_0.$$  
Then \eqref{then} combined with 
$$A\circ X_0=X_0\circ B_0$$
implies that  
$$ \delta \circ \t B_0 \circ \delta^{-1} = \mu_i \circ B_0$$  
for some $i$, $1 \leq i \leq n$.  
Consequently,  
\[
\t B_0 = \delta^{-1} \circ \mu_i \circ B_0 \circ \delta
\]  
for some $i$, $1 \leq i \leq n$. 
Thus,  if polynomials $A$, $B_0$, and $X_0$ form a primitive
solution of \eqref{1} with fixed $n\geq 2$ and fixed $\nu(0)$, then, up to
conjugacy, there are at most $n$ possibilities for $B_0$. 
 \qed

\bt \l{fin00}  For any $d \ge 2$, there exists a constant $r(d)$ such that for every  polynomial $R$ of degree $d$, there are polynomials
$B_1, \dots, B_r$  with $r \le r(d)$ satisfying
the following property: for a polynomial $B$, the relation
$B^{\circ l}\leq \delta\circ R^{\circ n},$ where $l,n\geq 1$ and $\delta\in \Sigma(R)$, implies that $B$ is conjugate to $\mu\circ B_i^{\circ s}$
for some $i$, $1\leq i\leq r$, some $\mu\in\Sigma(B_i)$,
and some $s\geq 1$.

\et  

\pr 
If $R$ is special, Theorem~\ref{spe} implies that the theorem holds for $r=1$ and  any polynomial $B_1$  that shares its Julia set with $R$, and whose degree $d_0$ is the minimal natural number such that $d = d_0^k$ for some $k \ge 1$.
Thus, it is enough to prove the theorem for non-special $R$.

By Theorem~\ref{t30}, the relation \be \l{reel} B^{\circ l}\leq \delta\circ R^{\circ n}\ee implies that there exists a polynomial
$\widehat{B}$ such that
\be
\widehat{B} \le R
\qquad\text{and}\qquad
J(\widehat{B}) = J(B).
\ee
On the other hand, by Theorem~\ref{ko1}, there exist $r(d)$ and  polynomials
$\widehat{B}_1, \widehat{B}_2, \dots, \widehat{B}_r$  with $r \le r(d)$ such that $\h B\leq  R$ implies that $\widehat{B}$ is conjugate to one of $\widehat{B}_1, \widehat{B}_2, \dots, \widehat{B}_r$.
 Thus, the relation \eqref{reel} 
implies that 
\[
J(B) = J\!\left(\alpha \circ \widehat{B}_i \circ \alpha^{-1}\right)
\]
for some $i,$ $1 \le i \le r$, and some polynomial $\alpha$ of degree one,
or, equivalently, that 
\be \l{byby} 
J\!\left(\alpha^{-1} \circ B \circ \alpha\right) = J(\widehat{B}_i).
\ee

Since $R$ is non-special, Theorem~\ref{spe} implies that $\h B$ and  
$\widehat{B}_1, \widehat{B}_2, \dots, \widehat{B}_r$ are also non-special.  
Therefore, by Theorem~\ref{sha}, there are polynomials  
$B_1, B_2, \dots, B_r$ such that the equality  
\[
J(P)=J(\widehat{B}_i), \qquad 1 \le i \le r,
\]  
for a polynomial $P$ implies  that 
\[
P=\delta \circ B_i^{\circ s}, \qquad 1 \le i \le r,
\]  
for some $\delta \in \Sigma(B_i)$ and some $s \ge 1$.  
By \eqref{byby}, this completes the proof.  
\qed

\bt \l{fin}  
For any $d \ge 2$, there exists a constant $r(d)$ such that for every
 polynomial $A$ of degree $d$, there are polynomials
$R_1, R_2, \dots, R_r$ of degree at most $d$ with $r \le r(d)$
satisfying the following property: for any polynomial $R$, the relation
$
A^{\circ k} \le  R,
$ 
where  $k\geq 1$,
implies that $R$ is conjugate to $\mu \circ R_i^{\circ n}$ for some $i$, $1 \le i \le r$, some $\mu \in \Sigma(R_i)$, and 
some  $n \ge 1$.

\et  

\pr If $A$ is special, Theorem~\ref{spe} implies that the conclusion of the theorem holds for $r=1$ and $R_1 = A$.  
Thus, we may assume that $A$ is not special.

The relation
\be \l{rle}
A^{\circ k}\leq R
\ee
implies by Theorem~\ref{t33} that there exist an integer $e\geq1$
and a polynomial $T$ such that
\be \l{ina0}
R^{\circ e}=T^{\circ ek}.
\ee
Hence,
\[
A^{\circ ke} \le T^{\circ ek},
\]
implying by Theorem~\ref{t30} that
there exists a polynomial \(\widehat{A}\) such that
$$\h A \le T \qquad\text{and}\qquad J(\h A)=J(A).$$
Furthermore, since
$$\deg \widehat{A} = \deg A,$$
it follows from Corollary~\ref{shaa} that
\[
\widehat{A} = \gamma \circ A
\]
for some \(\gamma \in \Sigma(A)\).

Since Theorem~\ref{cic} yields an explicit bound
$\vert \Sigma(A)\vert \leq d$, it follows from
Theorem~\ref{k02} that there exist $r(d)$ and polynomials
$T_1, T_2, \dots, T_r$ of degree $d$ with $r \le r(d)$
such that the condition
\be \l{tf0}
\gamma \circ A \leq T, \quad \gamma\in \Sigma(A),
\ee
implies that $T$ is conjugate to one of $T_1, T_2, \dots, T_r$.
Since \eqref{ina0} implies $J(R)=J(T)$, we conclude that 
the relation \eqref{rle} implies that 
\[
J(R) = J\!\left(\alpha \circ T_i \circ \alpha^{-1}\right)
\]
for some $i,$ $1 \le i \le r$, and some polynomial $\alpha$ of degree one,
or, equivalently, that 
\be \l{ybyb}
J\!\left(\alpha^{-1} \circ R \circ \alpha\right) = J({T}_i).
\ee 

Since $A$ is not special, Theorem~\ref{spe} implies that $T$ and $T_1, T_2, \dots, T_r$
are not special either. Therefore, by Theorem~\ref{sha}, there exist polynomials  
\(R_1, R_2, \dots, R_r\) of degree at most \(d\) such that the equality  
\[
J(P)=J(T_i), \quad 1 \le i \le r,
\]
for a polynomial \(P\) implies that  
\[
P=\delta \circ R_i^{\circ s}
\]
for some \(\delta \in \Sigma(R_i)\) and some \(s \ge 1\).  
By \eqref{ybyb}, this completes the proof. 
\qed

\vskip 0.2cm
\noindent{\it Proof of Theorem \ref{t2}.} Suppose first that $A$ is special, and let
$B\in\operatorname{Inter}(A)$. Then $B$ is also special,
although it need not have the same special type as $A$
(see \cite[Proposition~5.1]{gn1} or the end of
Section~\ref{24}).
Hence, $B$ is conjugate to $z^{d_0^b}$ or $\pm T_{d_0^b}$
for some integer $b\geq 1$,
where $d_0\geq 2$ is the smallest integer such that
$\deg A=d_0^a$ for some integer $a\geq 1$. Thus, the theorem holds in this case with
$B_1=z^{d_0}$ and $B_2=T_{d_0}$, since 
 $-z\in\Sigma(T_{d_0})$.

Assume now that $A$ is non-special, in which case
every $B\in\operatorname{Inter}(A)$ is also non-special. The condition $B \in \mathrm{Inter}(A)$ means that
\be \l{re1} 
A^{\circ k} \le R, \quad B^{\circ l} \le R
\ee
for some polynomial $R$ and integers $k, l \ge 1$. 

By Theorem~\ref{fin}, there exist  polynomials
$R_1, R_2, \dots, R_r$ of degree at most $d$ with $r \le r_1(d)$, 
 where $r_1(d)$ depends only on $d$,  
such that
the first relation in \eqref{re1} 
 implies that $R$ is conjugate to $\delta \circ R_i^{\circ n}$
for some $i$, $1 \le i \le r$, some $\delta \in \Sigma(R_i)$, and some  $n \ge 1$. On the other hand,  
Theorem~\ref{fin00} yields that  for each fixed $i,$ $1\leq i \leq r$,  there exist
polynomials 
$B_1,\dots,B_{r'}$ with $r' \le r_2(d)$, where $r_2(d)$ depends only on $d$,  such that the relation
$$B^{\circ l}  \le \delta \circ R_i^{\circ n},$$ where $\delta \in \Sigma(R_i)$ and  $n\geq 1,$ implies that $B$ is conjugate to $\mu\circ B_j^{\circ s}$
for some $j$, $1\leq j\leq r'$, some $\mu\in\Sigma(B_j)$,
and some $s\geq 1$.

Since, for any degree-one polynomial $\gamma$ and any polynomial $V$, the set of polynomials $U$ satisfying $U \le V$ coincides with the set of polynomials $U$ satisfying $$U \le \gamma \circ V \circ \gamma^{-1},$$ the theorem follows.  
\qed

\section*{acknowledgements}
The author is grateful to Geng-Rui Zhang for his comments on the paper.

\bibliographystyle{amsplain}

\end{document}